\documentclass[12pt,a4paper]{amsart}
\date{9 December 2010}

\usepackage{latexsym,amsmath,amsfonts,amscd,amssymb, mathrsfs, setspace, slashed, amsthm}
\usepackage[english]{babel}

\theoremstyle{plain}  
\newtheorem{theorem}{Theorem}[section]

\newtheorem*{theorem*}{Theorem}

\newtheorem{lemma}[theorem]{Lemma}

\theoremstyle{remark}

\newtheorem{remark}[theorem]{Remark}

\newtheorem*{claim*}{Claim}

\numberwithin{equation}{section}

\newcommand{\n}{\noindent}

\renewcommand{\leq}{\leqslant}

\newcommand{\lra}{\longrightarrow}
\newcommand{\R}{\mathbb{R}}

\renewcommand{\S}{\slashed{\mathcal{S}}}

\def\proof{\noindent\textit{Proof ---}\hspace*{0.2cm}}
\def\qed{\vspace*{-0.1cm} \hfill{$\square$}}
\def\disfrac#1#2{\displaystyle{\frac{#1}{#2}}}
\newcommand{\D}{\slashed{D}}

\def\End{\mathrm{End}}

\begin{document}
\title{A vanishing theorem in twisted de Rham cohomology}

\author[A. C. Ferreira]{Ana Cristina Ferreira}
\address{Departmento de Matem\'{a}tica \\
Universidade do Minho \\
Campus de Gualtar \\
4710-057 Braga \\
Portugal} \email{anaferreira@math.uminho.pt}

\thanks{Partially funded by FCT through the POPH-QREN scholarship program. Partial financial support provided by the Research Centre of Mathematics of the University of Minho through the FCT pluriannual funding program}

\begin{abstract}
We prove a vanishing theorem for the twisted de Rham cohomology of a
compact manifold.
\end{abstract}

\maketitle

\section{Introduction}

In this article, we show how to use connections with skew torsion to
identify the operator $(d+H)+(d+H)^*$, where $H$ is a three-form, with a cubic
Dirac operator. In the compact case, if $H$ is
closed, we prove a vanishing theorem for twisted de Rham cohomology
by means of a Lichnerowicz formula. As an application, we prove that
for a compact non-abelian Lie group the cohomology of the complex
defined by $d+H$, where $H$ is the three-form defined by the Lie
bracket, vanishes.

\bigskip
\section{The Dirac operator}

Let $(M,g)$ be a Riemannian manifold. Suppose that $\nabla$ is a connection on the tangent
bundle of $M$ and let $T$ be its (1,2) torsion tensor. If we
contract $T$ with the metric we get a (0,3) tensor which we will
still call the torsion of $\nabla$. If $T$ is a three-form then we
say that $\nabla$ is a connection with skew-symmetric torsion. Given
any three-form $H$ on $M$ then there exists a unique metric
connection with skew torsion $H$  defined explicitly by
$$g(\nabla_X Y, Z) = g(\nabla^g_X Y, Z)+\mbox{\small{$\frac{1}{2}$}} H(X,Y,Z)$$ where $\nabla^g$ is the Levi-Civita connection.

Fix a three-form $H$ and consider the one-parameter family of affine
connections
$$\nabla^s := \nabla^g + 2 s H$$ (Notice that if $s=\frac{1}{4}$ we
recover the connection with torsion $H$.)  If $M$ is spin, these connections lift to
the spin bundle $\S$ of $M$ as
$$\nabla^s_X (\varphi) := \nabla^g_X (\varphi) + s (i_X H)\varphi$$ where
$X$ is a vector field, $\varphi$ is a spinor field, and $i_X H$ is
acting by Clifford multiplication.

We may define the Dirac operator $\D$ on $\S$  with respect to
$\nabla$ by means of the following composition
$$\Gamma(M,\S)\lra \Gamma(M,T^*M\otimes\S)\lra \Gamma(M,TM\otimes\S)\lra \Gamma(M,\S)$$
where the first arrow is given by the connection, the second by the
metric and the third by the Clifford action. Suppose now that we have
a complex vector bundle $\mathcal{W}$, we can form the tensor
product $\S\otimes\mathcal{W}$, which is usually called a twisted
spinor bundle or a spinor bundle with values in $\mathcal{W}$. If
$\mathcal{W}$ is equipped with a Hermitian connection
$\nabla^\mathcal{W}$, we can consider the tensor product connection
$\nabla\otimes 1 + 1 \otimes \nabla^\mathcal{W}$, again denoted by
$\nabla$, on $\S\otimes\mathcal{W}$. We can define a Dirac
operator on this twisted spinor bundle associated with the
connection $\nabla$ by the same formula, where the action of the
tangent bundle by Clifford multiplication is only on the left
factor.

We will need to make use of a Lichnerowicz type formula for the
square of the Dirac operator. Such a formula first appeared in the
literature in \cite{Bismut}. See also \cite{AgricFried}.

\begin{theorem}\label{theo: Lich} {\bf [Bismut, \cite{Bismut}]}
The rough Laplacian $\Delta^s = {\nabla^s}^* \nabla^s$ and the square of the Dirac operator
$\mathrm{D}^{s/3}$ are related by
$$\left(\mathrm{D}^{s/3}\right)^2 = \Delta^s + F^\mathcal{W} +\disfrac{1}{4} \kappa + s dH - 2 s^2 ||H||^2,$$
where $\kappa$ is the Riemannian
scalar curvature and $F$ is the curvature of the twisting bundle
acting as $\sum_{i<j}F^\mathcal{W}(e_i,e_j)e_ie_j$ on
$\S\otimes\mathcal{W}$.
\end{theorem}

Notice that this formula relates the square of the Dirac operator
$\mathrm{D}^{s/3}$ and the Laplacian $\Delta^s$. The Dirac operator $\mathrm{D}^{1/3}$
is usually referred to as the cubic Dirac operator.

\bigskip

\section{Twisted cohomology}

Consider the spinor bundle with values in itself, that is,
$\S\otimes\S$. Recall that for this we do not need a global spin structure. We have, in even dimensions, the following chain of
isomorphisms
$$\S\otimes\S \simeq \S^*\otimes\S \simeq \End (\S) \simeq \mathrm{Cl}
\simeq \Lambda$$ where $\mathrm{Cl}$ denotes the Clifford bundle and $\Lambda$ the bundle of exterior forms.

If we take the induced Levi-Civita connection $\nabla^g$ on both
factors of $\S\otimes\S$ and consider the tensor product connection
 $\nabla^g\otimes 1 + 1\otimes \nabla^g$ we obtain the
induced Levi-Civita connection, again denoted by $\nabla^g$, on
$\Lambda$. If we consider the associated Dirac operator $\mathrm{D}^g$ on
$\S\otimes \S$ we get a familiar operator on $\Lambda$. In fact,
$$\mathrm{D}^g = d+d^\ast$$ where $d$ is the exterior differential and $d^\ast$ is its formal adjoint, \cite{Lawson}.

The same fact can be claimed for an odd-dimensional manifold.
Consider the inclusion $M\hookrightarrow \R\times M$, $\S^+$ and $\S^-$ the half spinor bundles of $\R\times M$. The Clifford action by $e_0$, where $e_0$ is a unit vector field of $\R$, gives an isomorphism between $\S^+$ and $\S^-$ and thus we can regard $\S^+\simeq \S^-$ as the spinor bundle of $M$.
Under this identification, the
Dirac operator associated to the Levi-Civita connection becomes
$$\S^+\stackrel{\mathrm{D}^g}{\lra}\S^- \stackrel{e_0}{\lra} \S^+$$
where $e_0$ denotes multiplication by $e_0$. Consider also the Levi-Civita connection on $\S$ and
the twisted Dirac operator
$$\S^+\otimes \S \stackrel{\mathrm{D}^g}{\lra} \S^-\otimes\S \stackrel{e_0}{\lra} \S^+\otimes \S.$$
Notice that the exterior bundle of $M$ is $\Lambda \simeq
\mathrm{Cl}\simeq\S^+\otimes\S$, and so the twisted Dirac operator
above is, in terms of differential forms, the restriction of the
operator $d+d^*$ on $\mathbb{R}\times M$ to forms that are
independent of the coordinate $t$ of $\mathbb{R}$, and can therefore
be seen as $d+d^*$ on $M$.

We may now ask ourselves what happens if we introduce connections
with skew torsion in this setting.

\begin{theorem}\label{theo: twisted}
Let $H$ be a three-form, and suppose that the left and right spinor
factors are, respectively, equipped with the connections $\nabla^g +
\frac{1}{12}H$ and $\nabla^g - \frac{1}{4} H$. Consider the tensor
product of these two connections on $\S\otimes\S$. The corresponding
Dirac operator on $\Lambda$ is given by
$$\mathrm{D}= (d+H)+(d+H)^*$$
where $H$ is acting by exterior multiplication and $(d+H)^*$ is the
formal adjoint of $d+H$ with respect to the metric, namely, $d^* + (-1)^{n(p+1)}
* H *$ on $\Lambda^p $.
\end{theorem}

\proof Let us consider first an even dimensional manifold. Take a $p$-form $\theta$ and identify it with $\varphi= \sum_r \varphi_r^+\otimes\varphi_r^- \in
\Gamma(M,\S\otimes\S).$ Then the Clifford left and right actions
of a vector field $e$ are given, respectively, by
$$\begin{array}{lclclcl} e \varphi & = & \sum_r e\varphi_r^+ \otimes \varphi_r^- & = & e\wedge \theta - e \lrcorner \theta\\ \varphi e & = & \sum_r \varphi_r^+ \otimes e \varphi_r^- & = & (-1)^p (e\wedge \theta + e \lrcorner \theta)\end{array}$$

Using the summation convention, we have
$$\begin{array}{lcl} \mathrm{D}(\varphi) & = & e_i \nabla^g_{e_i} \varphi_r^+ \otimes \varphi_r^-
+  e_i \varphi_1 \otimes \nabla^g_{e_i} \varphi_2 + \medskip \\ & &
\frac{1}{12} e_i (e_i \lrcorner H) \varphi_r^+ \otimes \varphi_r^-  -
\frac{1}{4} e_i
\varphi_r^+ \otimes (e_i \lrcorner H) \varphi_r^- \medskip \\
& = & e_i\nabla^g_{e_i}(\varphi) +  \frac{1}{12} e_i(e_i \lrcorner
H) \varphi +\frac{1}{4} e_ i \varphi (e_i\lrcorner H).\end{array} $$
Since $\mathrm{D}^g(\varphi) = e_i \nabla^g_{e_i} (\varphi)$ corresponds to $(d+d^*)\theta$, it remains to see that
$\frac{1}{12} e_i(e_i\lrcorner H)\varphi + \frac{1}{4} e_i\varphi (e_i\lrcorner \varphi)$ can be identified with $(H+H^*)\theta$.

\smallskip

Write $H =  H_{abc} e_a \wedge e_b \wedge e_c$ and
observe that $$H_{abc} e_a\wedge e_b \wedge e_c \wedge \theta + H_{abc} e_c\lrcorner (e_b \lrcorner (e_a \lrcorner \theta))
$$ is the same as $(H+H^*)\theta$ since the formal adjoint of exterior multiplication is interior multiplication.
It is simple to see that $e_i (e_i\lrcorner H)\varphi = 3 H \varphi$ and that the action of $H$ is given by
$$H_{abc} (e_a \wedge e_b \wedge e_c \wedge \theta + e_a \wedge e_b \wedge (e_c \lrcorner \theta) + e_a \wedge (e_b
\lrcorner (e_c \lrcorner \theta) + \dots ) $$
and that $e_i\varphi (e_i\lrcorner H)\theta$ is such that when we add
$$\frac{1}{12}e_i(e_i\lrcorner H)\theta = \frac{1}{4}H\theta$$ and $$\frac{1}{4}e_i \theta (e_i\lrcorner H)$$ the mixed terms cancel and it amounts to
$$\frac{1}{4} H_{abc} ( e_a\wedge e_b \wedge e_c \wedge \theta +  e_c\lrcorner (e_b \lrcorner (e_a \lrcorner \theta))$$ plus $$\frac{3}{4} H_{abc} (e_a\wedge e_b \wedge e_c \wedge \theta + e_c\lrcorner (e_b \lrcorner (e_a \lrcorner \theta))$$ which is then  $(H+H^*)\theta$. The proof in the odd-dimensional case is perfectly analogous.

\qed

\begin{remark}
Notice that these are lifts of the metric connections on the tangent
bundle with torsion $\frac{1}{3}H$ and $-H$. It is interesting to
observe that these weights $\frac{1}{3}$ and $-1$ also appear in
Bismut's proof of the local index theorem for non-K\"{a}hler
manifolds, \cite{Bismut}.
\end{remark}

Suppose now that $H$ is a closed three-form. In \cite{AtiyahSegal},
Atiyah and Segal defined the concept of twisted de Rham cohomology.
On the de Rham complex of differential forms $\Omega$ we can define
the operator $d+H$. Note that $$(d+H)^2 = d^2 + d H + H d + H^2 =
0$$ since $H$ is closed and of odd degree. The operator $d+H$ does
not preserve form degrees but preserves the $\mathbb{Z}_2$-grading.
We then have a 2-step chain complex and the cohomology of this
complex is then the twisted de Rham cohomology.

The twisted de Rham complex is an elliptic complex so, on a compact manifold, Hodge theory applies.
If $H^+$ and $H^-$ are the cohomology groups then
$$H^\pm \simeq
\{\theta \in \Omega^\pm: ~(d+H)\theta = 0 \quad \mbox{and} \quad (d+H)^\ast \theta = 0\}$$ or,
in other words, each cohomology class has a unique representative in the kernel of $\mathrm{D}^2$ where $$\mathrm{D} = (d+H)+(d+H)^\ast.$$

\bigskip

\section{A vanishing theorem}

We can use the Lichnerowicz formula of theorem \ref{theo: Lich} and
also theorem \ref{theo: twisted} to prove the following

\begin{theorem}\label{theo: vanishing}
Let $M$ be a compact spin manifold and let $H$ be a closed
three-form. Consider the Dirac operator $\mathrm{D}^{1/12}$ on $\S\otimes \S$
associated with the connection $$\nabla = \nabla^{1/12} \otimes 1 +
1 \otimes \nabla^{-1/4},$$  let $F^{-1/4}$ be the curvature of
$\nabla^{-1/4}$ on $\S$ and $\kappa$ the Riemannian scalar curvature
of $M$. If
$$F^{-1/4} + \frac{1}{4}\kappa - \frac{1}{8} \Vert H \Vert^2$$ acts as a
positive endomorphism then the twisted de Rham cohomology for $d+H$
vanishes.
\end{theorem}

\proof We start by observing that we need only to prove that the kernel of the operator $\mathrm{D}^{1/12}$ is zero.
Consider $\psi$ a smooth section of $\S\otimes\S$. Since $dH=0$,
the Lichnerowicz formulas gives
$$\left(\mathrm{D}^{1/12}\right)^2\psi = \Delta^{1/4}\psi + F^{-1/4}\psi+ \frac{1}{4} \kappa \psi -
\frac{1}{8} \Vert H \Vert^2 \psi. $$ Now take the inner
product of this with $\psi$. Since the Dirac
operator is self-adjoint and the Laplacian $\Delta$ is given by
$\nabla^*\nabla$, we get
\small{$$\int_M  \Vert \mathrm{D}^{1/12}\psi \Vert^2 ~ \mathrm{dVol}  = \int_ M  \Vert \nabla^{1/4}\psi \Vert^2 + (F^{-1/4}\psi, \psi)+ \frac{1}{4} \kappa \Vert \psi \Vert^2  -
\frac{1}{8} \Vert H \Vert^2 \Vert \psi \Vert^2 ~\mathrm{dVol}.$$}
Using the hypothesis that
$$F^{-1/4} + \frac{1}{4}\kappa - \frac{1}{8} \Vert H \Vert^2$$ is a
positive endomorphism we conclude that $\mathrm{D}^{1/12}\psi = 0$ if and only if
$\psi = 0$.

\qed

\bigskip

\section{An example}

Let $G$ be a compact, non-abelian Lie group equipped with a
bi-invariant metric. Consider the one-parameter family of
connections $\nabla^t_X(Y) = t [X,Y]$. Given $t$, the torsion of
$\nabla^t$ is $(2t-1) [X,Y]$. Notice that since the metric is
ad-invariant, it means that these are metric connections and also
that their torsion is skew-symmetric. Note also that if
$t=\frac{1}{2}$ we get the Levi-Civita connection, since the torsion
vanishes. The curvature of $\nabla^t$ is given by
$$R^{\nabla^t}(X,Y)Z = t^2[X,[Y,Z]]-t^2[Y,[X,Z]]-t[[X,Y],Z] = (t^2-t)[[X,Y],Z],$$
by means of the Jacobi identity. For $t=0$ and $t=1$, we get two
flat connections. These correspond, respectively, to the left and
right invariant trivialization of the tangent bundle, \cite{Nomizu}.

Let us write the above one-parameter family of connections as
$$\nabla^{2s}_X(Y) = \nabla^g_X(Y) + 2s [X,Y].$$
Notice that the Levi-Civita connection corresponds now to the
parameters $s=0$ while the two flat connections correspond to $s=\pm
\frac{1}{4}$.

Consider the lift of these connections to the spinor bundle $\S$ of
$G$. Take the connection $\nabla^{1/12}\otimes 1 + 1\otimes
\nabla^{-1/4}$ on $\Gamma(M,\S\otimes \S)$. We know from theorem
\ref{theo: twisted} that the Dirac operator $\mathrm{D}^{1/12}$ then
corresponds to $(d+H)+(d+H)^*$ on $\Lambda G$, where $H$ is given by
$H(X,Y,Z)=([X,Y],Z).$ Note that $H$, being a bi-invariant form, is
closed.

We need the following auxiliary lemma, which can be proved by direct
computation.

\begin{lemma}\label{Lemma: Lie-scalar}
Let $G$ be a non-abelian Lie group equipped with a bi-invariant
metric, then the scalar curvature $\kappa$ of $G$ is given by
$$\kappa = \frac{1}{4}
\sum_{ij} \Vert [e_i,e_j] \Vert ^2 $$ where $\{e_i \}$ is an orthonormal basis of the Lie algebra of $G$.
\end{lemma}

\begin{theorem}\label{Theo: Lie}
Let $G$ be a compact, non-abelian Lie group equipped with a
bi-invariant metric and let $H(X,Y,Z) = ([X,Y],Z)$ be the associated
bi-invariant three-form. Then the twisted de Rham cohomology of
$d+H$ vanishes.
\end{theorem}

\n \proof Since $F^{-1/4}= 0$, by means of theorem \ref{theo:
vanishing} we only need to show that the constant $\rho =
\frac{1}{4} \kappa - \frac{1}{8} \Vert H \Vert ^2$ is positive. We
have already computed $\kappa$ in lemma \ref{Lemma: Lie-scalar}, so
if we take the same orthonormal basis we get that
$$\Vert H \Vert ^2 = \frac{1}{6} \sum_{ijk}|([e_i,e_j],e_k)|^2,$$
and using the Cauchy-Schwarz inequality
$$\Vert H \Vert ^2 \leq \frac{1}{6} \sum_{ijk} \Vert [e_i,e_j] \Vert^2 \Vert e_k \Vert ^2 = \frac{1}{6} \sum_{ij}  \Vert [e_i,e_j] \Vert^2 $$
So $\rho > \left( \frac{1}{16}-\frac{1}{48} \right) \sum_{ij} \Vert
[e_i,e_j] \Vert^2  > 0 $.

\qed

\bigskip
\begin{remark}
To see this result for connected, compact, simple groups in a
different way, note that it is well known that by averaging, each
cohomology class of $G$ can be represented by a bi-invariant form.
The de Rham cohomology ring $H^*(G)$ is an exterior algebra
(more precisely $H^*(G)$ is an exterior algebra on generators in
degree $2d_i-1$, where each $d_i$ is the degree of generators of
invariant polynomials on the Lie algebra of $G$). The Killing form gives
$H^3(G) = \mathbb{R}$. Consider now the twisted de Rham operator
$d+H$. Since $H$ is bi-invariant, the twisted cohomology classes can also be represented by bi-invariant forms.
Since bi-invariant forms are closed, $(d+H)\alpha = H\wedge \alpha$. So if $H\wedge \alpha = 0$, since $H$ is a generator, then $H\wedge \alpha =0$ implies that $\alpha = H\wedge \beta$ for some $\beta$.
Therefore, the twisted cohomology vanishes.
\end{remark}

\bigskip

\noindent   {\bf Acknowledgements:}   I would like to thank Nigel Hitchin for pointing me to this topic, and
for the many helpful conversations that have ensued.

\end{document}